\documentclass[11pt,leqno]{article}
\usepackage{amsmath}
\usepackage{amsthm}

\numberwithin{equation}{section}
\let\sect=\section

\newtheorem{theorem}{Theorem}[section]
\newtheorem{corollary}[theorem]{Corollary}
\newtheorem{lemma}[theorem]{Lemma}
\newtheorem{proposition}[theorem]{Proposition}
\newtheorem{claim}[theorem]{Claim}
\newtheorem{example}[theorem]{\sl Example}

\theoremstyle{definition}
\newtheorem{remark}[theorem]{Remark}

\newcommand{\EE}{{\bf  E}}

\newcommand{\RR}{{\bf  R}}

\newcommand{\Var}{{\bf Var}}

\newcommand{\Lc}{{\cal L}}

\newcommand{\Leq}{{\,\stackrel{\Lc}{=}\,}}

\newcommand{\begp}{\begin{proposition}}
\newcommand{\enp}{\end{proposition}}
\newcommand{\begt}{\begin{theorem}}
\newcommand{\ent}{\end{theorem}}
\newcommand{\begl}{\begin{lemma}}
\newcommand{\enl}{\end{lemma}}
\newcommand{\begc}{\begin{corollary}}
\newcommand{\enc}{\end{corollary}}
\newcommand{\begcl}{\begin{claim}}
\newcommand{\encl}{\end{claim}}
\newcommand{\begr}{\begin{remark}}
\newcommand{\enr}{\end{remark}}
\newcommand{\begal}{\begin{algorithm}}
\newcommand{\enal}{\end{algorithm}}
\newcommand{\begd}{\begin{definition}}
\newcommand{\enf}{\end{definition}}
\newcommand{\begx}{\begin{example}}
\newcommand{\enx}{\end{example}}
\newcommand{\bega}{\begin{array}}
\newcommand{\ena}{\end{array}}

\def\rompar(#1){\textup(#1\textup)}    
\newcommand\xfrac[2]{#1/#2}
\newcommand\eps{\varepsilon}
\newcommand\cS{\mathcal S}

\newcommand{\refT}[1]{Theorem~\ref{#1}}
\newcommand{\refC}[1]{Corollary~\ref{#1}}
\newcommand{\refL}[1]{Lemma~\ref{#1}}
\newcommand{\refR}[1]{Remark~\ref{#1}}

\newcommand\eg{e.g.\spacefactor=1000}
\newcommand\cf{{cf.}\spacefactor=1000}

\renewcommand\Re{\operatorname{Re}}

\newcommand\nopf{\qed}   
\newcommand\noqed{\renewcommand{\qed}{}} 

\newcommand\Roesler{R\"{o}sler}

\begin{document}

\setcounter{page}{0}
\thispagestyle{empty}

\begin{center}
{\Large \bf Smoothness and Decay Properties of \\ the Limiting Quicksort Density Function
\\ } 
\normalsize
 
\vspace{4ex}
{\sc James Allen Fill\footnotemark} \\
\vspace{.1in}
Department of Mathematical Sciences \\ 
\vspace{.1in}
The Johns Hopkins University \\
\vspace{.1in}
{\tt jimfill@jhu.edu} and {\tt http://www.mts.jhu.edu/\~{}fill/} \\
\vspace{.2in}
{\sc and} \\
{\sc Svante Janson}\\ 
\vspace{.1in}
Department of Mathematics \\ 
\vspace{.1in}
Uppsala University \\
\vspace{.1in}
{\tt svante.janson@math.uu.se} and {\tt http://www.math.uu.se/\~{}svante/} \\
\end{center}
\vspace{3ex}
 
\begin{center}
{\sl ABSTRACT} \\
\end{center}
\vspace{.05in}

\begin{small}
Using Fourier analysis, we prove that the limiting distribution of the standardized
random number of comparisons used by {\tt Quicksort} to sort an array of~$n$ numbers
has an everywhere positive and infinitely differentiable density~$f$, and that each
derivative~$f^{(k)}$ enjoys superpolynomial decay at~$\pm \infty$.  In particular,
each~$f^{(k)}$ is bounded. Our method is sufficiently computational to prove, for example,
that~$f$ is bounded by~$16$.
\smallskip
\par\noindent
{\em AMS\/} 2000 {\em subject classifications.\/}  Primary 68W40;
secondary 68P10, 60E05, 60E10.
\medskip
\par\noindent
{\em Key words and phrases.\/} {\tt Quicksort}, density, characteristic function,
sorting algorithm, Fourier analysis, rapidly decreasing~$C^{\infty}$ function,
tempered distribution, integral equation.
\medskip
\par\noindent
\emph{Date.\/} March 31, 2000.
\end{small}

\footnotetext[1]{Research supported by NSF grant DMS--9803780,
and by the Acheson J.~Duncan Fund for the Advancement of Research in Statistics.}

\newpage
\addtolength{\topmargin}{+0.5in}

\sect{Introduction and summary}

The {\tt Quicksort} algorithm of Hoare~\cite{Hoare} is ``one of the fastest, the
best-known, the most generalized, the most completely analyzed, and the most
widely used algorithms for sorting an array of numbers''~\cite{ES}.  {\tt Quicksort}
is the standard sorting procedure in {\tt Unix} systems, and Philippe Flajolet,
a leader in the field of analysis of algorithms, has noted that it is among ``some
of the most basic algorithms---the ones that do deserve deep
investigation''~\cite{Flaj}.
Our goal in this introductory section is to review briefly some of what is known about
the analysis of {\tt Quicksort} and to summarize how this paper advances that analysis.

The {\tt Quicksort} algorithm for sorting an array of~$n$ numbers is extremely simple
to describe.  If $n = 0$ or $n = 1$, there is nothing to do.  If $n \geq 2$, pick a number
uniformly at random from the given array.  Compare the other numbers to it to partition the
remaining numbers into two subarrays.  Then recursively invoke {\tt Quicksort} on each of
the two subarrays.

Let~$X_n$ denote the (random) number of comparisons required (so that $X_0 = 0$).  Then~$X_n$
satisfies the distributional recurrence relation
$$
X_n \Leq X_{U_n - 1} + X^*_{n - U_n} + n - 1,\qquad n \geq 1,
$$
where~$\Leq$ denotes equality in law (i.e.,\ in distribution), and where, on the right,
$U_n$ is distributed uniformly on the set $\{1, \ldots, n\}$,
$X_j^* \Leq X_j$, 
and
$$
U_n;\ X_0, \ldots, X_{n - 1};\ X^*_0, \ldots, X^*_{n - 1}
$$
are all independent.

As is well known and quite easily established, for $n \geq 0$ we have
$$
\mu_n := \EE X_n = 2 (n + 1) H_n - 4 n \sim 2 n \ln n,
$$
where $H_n := \sum_{k = 1}^n k^{-1}$ is the $n$th harmonic number and~$\sim$ denotes
asymptotic equivalence.  It is also routine to compute explicitly the standard deviation
of~$X_n$ (see Exercise 6.2.2-8 in~\cite{Knuth3}), which turns out to be
$\sim n \sqrt{7 - \frac{2}{3} \pi^2}$.

Consider the standardized variate
$$
Y_n := (X_n - \mu_n) / n,\qquad n \geq 1.
$$
R\'{e}gnier~\cite{Reg} showed using martingale arguments that $Y_n \to Y$
in distribution, with~$Y$ satisfying the distributional identity
\begin{equation}
\label{fix}
Y \Leq U Y + (1 - U) Z + g(U) =: h_{Y, Z}(U),
\end{equation}
where
\begin{equation}\label{gu}
g(u) := 2 u \ln u + 2 (1 - u) \ln (1 - u) + 1,
\end{equation}
and where, on the right of~$\Leq$ in~\eqref{fix},
$U$, $Y$, and~$Z$ are independent, with $Z \Leq Y$  and $U \sim
\mbox{unif}(0, 1)$.  R\"{o}sler~\cite{Roesler} showed that~\eqref{fix}
characterizes the limiting law~$\Lc(Y)$, in the precise sense that
$F := \Lc(Y)$ is the \emph{unique} fixed point of the operator
$$
G = \Lc(V) \mapsto SG := \Lc(U V + (1 - U) V^* + g(U))
$$
(in what should now be obvious notation) subject to
$$
\EE V = 0,\qquad \Var V < \infty.
$$

Thus it is clear that fundamental (asymptotic) probabilistic understanding of
{\tt Quick\-sort}'s behavior relies on fundamental understanding of the limiting
distribution~$F$.  In this regard, R\"{o}sler~\cite{Roesler} showed that
\begin{equation}
\label{finitemgf}
\mbox{the moment generating function (mgf) of~$Y$ is everywhere finite,}
\end{equation}
and Hennequin~\cite{H89}~\cite{H91} and R\"{o}sler showed how all the moments
of~$Y$ can be pumped out one at a time, though there is no known expression for the mgf
nor for the general $p$th moment in terms of~$p$.
Tan and Hadjicostas~\cite{TanH} proved that~$F$ has a density~$f$
which is almost everywhere positive, but their proof does not even show
whether~$f$ is continuous.

The main goal of this paper is to prove that~$F$ has a density~$f$ which is
infinitely differentiable, and that each derivative~$f^{(k)}(y)$ decays as $y \to \pm \infty$
more rapidly than any power of~$|y|^{-1}$: this is our main \refT{T:smooth}.
In particular, it follows that each~$f^{(k)}$ is bounded (cf.~\refT{T:maxf}).

Our main tool will be Fourier analysis.  We begin in Section~\ref{chf}
by showing (see \refT{T:chis}) that the characteristic function~$\phi$ for~$F$ has rapidly
decaying derivatives of every order.  Standard arguments reviewed briefly at the outset of
Section~\ref{density} then immediately carry this result over from~$\phi$ to~$f$.
Finally, in Section~4 we will use the boundedness and continuity of~$f$ to establish
an integral equation for~$f$ (\refT{T:inteq}).
As a corollary, $f$ is everywhere positive (\refC{C:positive}).

\begin{remark}
(a) Our method is sufficiently computational that we will prove, for example,
that~$f$ is bounded by~$16$.  This is not sharp numerically, as Figure~4 of~\cite{TanH}
strongly suggests that the maximum value of~$f$ is about $2 / 3$.  However, in future
work we will rigorously justify (and discuss how to obtain bounds on the error in)
the numerical computations used to obtain that figure, and the rather crude bounds
on~$f$ and its derivatives obtained in the present paper are needed as a starting point
for that more refined work.

(b) Very little is known rigorously about~$f$.  For example, the figure discussed in~(a)
indicates that~$f$ is unimodal.  Can this be proved?  Is~$f$ in fact \emph{strongly}
unimodal (i.e.,\ log-concave)?  What can one say about changes of signs 
for the derivatives of~$f$?

(c) Knessl and Szpankowski~\cite{KnSz} purport to prove very sharp estimates
of the rates of decay of $f(y)$ as $y \to - \infty$ and as $y \to \infty$.
Roughly put, they assert that the left tail of~$f$ decays doubly exponentially
(like the tail of an extreme-value density) and that the right tail decays
exponentially.  But their ``results'' rely on several unproven assumptions.
Among these, for example, is their assumption~(59) that
$$
\EE e^{- \lambda Y} \sim \exp( \alpha \lambda \ln \lambda + \beta \lambda + \gamma
\ln \lambda + \delta)\qquad \mbox{as $\lambda \to \infty$}
$$
for some constants $\alpha (> 0), \beta, \gamma, \delta$.
[Having assumed this, they derive the values of $\alpha$, $\gamma$, and~$\delta$ exactly,
and the value of~$\beta$ numerically.]
\end{remark}

\sect{Bounds on the limiting Quicksort characteristic function} \label{chf} 

We will in this section prove the following result on superpolynomial
decay of the characteristic function of the limit variable~$Y$.
\begin{theorem}\label{T:phi}
For every real $p \geq 0$ there is a smallest constant $0 < c_p < \infty$ such that
the characteristic function
$\phi(t) :\equiv \EE e^{itY}$ satisfies
\begin{equation}\label{phip}
|\phi(t)|\le c_p|t|^{-p}\mbox{\rm \ \ for all $t \in \RR$.}
\end{equation}
These best possible constants $c_p$ satisfy 
$c_0=1$,
$c_{1/2}\le2$,
$c_{3/4}\le\sqrt{8\pi}$, 
$c_1\le4\pi$, 
$c_{3/2}<187$, $c_{5/2}<103215$,
$c_{7/2}<197102280$,
and the relations
\begin{align}
c_{p_1}^{1/p_1}&\le c_{p_2}^{1/p_2},\qquad  0<p_1 \le p_2;
 \label{cpmon}\\
c_{p+1}&\le 2^{p+1} c_p^{1+(1/p)} p/(p-1),
 \qquad p>1;
 \label{cp+1}\\
c_p&\le 2^{p^2+6p},\qquad p>0.
 \label{cpbound}
\end{align}
\end{theorem}
[The numerical bounds are not sharp (except in the trivial case of~$c_0$); 
they are the best that we can get without too much work, but we expect that
substantial improvements are possible.]

\begin{proof}
The basic approach is to use the fundamental relation~\eqref{fix}.
We will first show,
using a method of van der Corput  \cite{vdC,Monty}, 
that the characteristic function of $h_{y,z}(U)$ is
bounded by $2|t|^{-1/2}$ for each $y, z$. Mixing, this yields
\refT{T:phi} for $p= 1 / 2$. Then we will use another consequence of 
\eqref{fix}, namely, the functional equation
\begin{equation}
\label{phiinteq}
\phi(t) = \int^1_{u = 0}\!\phi(u t)\,\phi((1 - u) t)\,e^{i t g(u)}\,du,\ \ t \in \RR,
\end{equation}
or rather its consequence
\begin{equation}
\label{phirecbd}
|\phi(t)| \leq \int^1_{u = 0}\!|\phi(u t)|\,\,|\phi((1 - u) t)|\,\,du,
\end{equation}
and obtain successive improvements in the exponent~$p$.

We give the details as a series of lemmas,
beginning with a standard calculus estimate \cite{Monty}.
Note that it suffices to consider $t > 0$ in the proofs because
$\phi(-t) = \overline{\phi(t)}$ and thus $|\phi(-t)|=|\phi(t)|$.
Note also that the best constants satisfy
$c_p = \sup_{t>0} t^p |\phi(t)|$ 
(although we do not know in advance of proving Theorem~\ref{T:phi} that these are
finite), and thus
$c_p^{1/p}=\sup_{t>0} t|\phi(t)|^{1/p}$,
which clearly satisfies~\eqref{cpmon} because $|\phi(t)|\le1$.

\begl
\label{calclem}
Suppose that a function~$h$ is twice continuously differentiable on an
open interval $(a, b)$ with
$$
h'(x) \geq c > 0\mbox{\rm \ \ \ and\ \ \ }h''(x) \geq 0
\mbox{\rm \ \ \ for $x \in (a, b)$.}
$$
Then
$$
\left| \int^b_{x = a}\!e^{i t h(x)}\,dx \right| \leq \frac{2}{c t}
\mbox{\rm \ \ for all $t > 0$.} 
$$
\enl

\begin{proof}
By considering subintervals $(a + \eps,b - \eps)$ and letting $\eps \to 0$,
we may without loss of generality assume that~$h$ is defined and 
twice differentiable at the endpoints, too.
Then, using integration by parts, we calculate
\begin{align*}
\int^b_{x = a}\!e^{i t h(x)}\,dx
 &= \frac{1}{i t} \int^b_{x = a}\!\left[ \frac{d}{dx} e^{i t h(x)}
 \right]\,\frac{dx}{h'(x)} 
      \\
 &= \frac{1}{i t} \left\{ \left. \frac{e^{i t h(x)}}{h'(x)}
 \right|^b_{x = a} 
- \int^b_{x =  a}\!e^{i t h(x)}\,d\left( \frac{1}{h'(x)} \right) \right\}.
\end{align*}
So
\begin{equation*}
\begin{split}
\left| \int^b_{x = a}\!e^{i t h(x)}\,dx \right|
 &\leq \frac{1}{t} \left\{ \left( \frac{1}{h'(b)} + \frac{1}{h'(a)} \right) 
  + \int^b_{x = a}\! \left| d \left( \frac{1}{h'(x)} 
  \right) \right|\,dx \right\} \\
 &  =  \frac{1}{t} \left\{ \left( \frac{1}{h'(b)} + \frac{1}{h'(a)}
   \right) + \int^b_{x = a}\! \left[ - d \left( \frac{1}{h'(x)} 
   \right) \right]\,dx \right\} \\
 &  =  \frac{1}{t} \left\{ \left( \frac{1}{h'(b)} + \frac{1}{h'(a)}
 \right) + \left(  \frac{1}{h'(a)} - \frac{1}{h'(b)} \right) \right\} \\
 &  =  \frac{2}{t h'(a)} \leq \frac{2}{c t}.
\hskip 75mm\qedsymbol\hskip-75mm
\end{split}
\end{equation*}
\noqed
\end{proof}

\begin{lemma}\label{L:hyz}
For any real numbers~$y$ and~$z$, the random variable~$h_{y,z}(U)$
defined by \eqref{fix} satisfies
\begin{equation*}
|\EE e^{ith_{y,z}(U)}|\le 2|t|^{-1/2}.
\end{equation*}
\end{lemma}

\begin{proof}
We will apply Lemma~\ref{calclem}, taking~$h$ to be~$h_{y, z}$.  Observe that
$$
h''_{y, z}(u) = 2 \left( \frac{1}{u} + \frac{1}{1 - u} \right) = \frac{2}{u (1 -
u)} \geq 8\mbox{\ \ for $u \in (0, 1)$}
$$
and that
$$
h'_{y, z}(u) = 0\mbox{\ \ \ if and only if\ \ \ }u = \alpha_{y, z} := \frac{1}{1 + \exp \left(
\frac{1}{2} (y - z) \right)} \in (0, 1).
$$
Let $t > 0$ and $\gamma > 0$.
If in Lemma~\ref{calclem} we take 
$a := \alpha_{y, z} + \gamma t^{-1/2}$ and $b := 1$, 
and assume that $a < b$, then note
$$
h'(u) = h'_{y, z}(u) = \int^u_{x = \alpha_{y, z}} h''_{y, z}(x)\,dx
\geq 8 (u - \alpha_{y,z}) 
\geq 8 \gamma t^{-1/2}\mbox{\ \ for all $u \in (a, b)$.}
$$
So, by Lemma~\ref{calclem},
$$
\left| \int^1_{u = \alpha_{y, z} + \gamma t^{-1/2}}\!e^{i t h_{y,
z}(u)}\,du \right| \leq 
\frac{2}{t} [8 \gamma t^{-1/2}]^{-1} = \frac{1}{4 \gamma} t^{-1/2}.
$$
Trivially,
$$
\left| \int^{\alpha_{y, z} + \gamma t^{-1/2}}_{u = \alpha_{y,
z}}\!e^{i t h_{y, z}(u)}\,du 
\right| \leq \gamma t^{-1/2},
$$
so we can conclude 
$$
\left| \int^{1}_{u = \alpha_{y, z}}\!e^{i t h_{y, z}(u)}\,du\right|
 \leq [ (4 \gamma)^{-1} +  \gamma] t^{-1/2}.
$$
This result is trivially also true when 
$a = \alpha_{y, z} + \gamma t^{-1/2} \ge b = 1,$
so it holds for all $t, \gamma > 0$.
The optimal choice of~$\gamma$ here is 
$\gamma_{\mbox{\scriptsize opt}} = 1 / 2$, 
which yields
$$
\left| \int^{1}_{u = \alpha_{y, z}}\!e^{i t h_{y, z}(u)}\,du\right|
 \leq  t^{-1/2}\qquad\text{for all }t>0.
$$
Similarly, for example by considering $u \mapsto h(1-u)$,
$$
\left| \int_{0}^{\alpha_{y, z}}\!e^{i t h_{y, z}(u)}\,du \right|
 \leq  t^{-1/2}\qquad\text{for all }t>0,
$$
and we conclude that the lemma holds for all $t>0$, and thus for all
real $t$.
\end{proof}

\begin{lemma}\label{L:1/2}
For any real $t$,
$|\phi(t)|\le 2|t|^{-1/2}$.
\end{lemma}

\begin{proof}
\refL{L:hyz} shows that 
$$
\Bigl|\EE\Bigl(e^{ith_{Y,Z}(U)}\Bigm| Y,Z\Bigr)\Bigr|
\le 2 |t|^{-1/2}
$$
and thus
$$
|\phi(t)|=
\Bigl|\EE e^{ith_{Y,Z}(U)}\Bigr|
\le \EE\Bigl|\EE\Bigl(e^{ith_{Y,Z}(U)}\Bigm| Y,Z\Bigr)\Bigr|
\le 2 |t|^{-1/2}.
\eqno\qedsymbol
$$
\noqed
\end{proof}

The preceding lemma is the case $p = 1/2$ of~\refT{T:phi}.
We now improve the exponent.

\begin{lemma}\label{L:p01}
Let $0<p<1$.  Then
$$
c_{2 p} \le \frac{\bigl[\Gamma(1-p)\bigr]^2}{\Gamma(2-2p)} c_p^2.
$$
\end{lemma}

\begin{proof}
By~\eqref{phirecbd} and the definition of~$c_p$,
\begin{equation*}
|\phi(t)|
\le
\int_{u = 0}^1 c_p^2|ut|^{-p}|(1-u)t|^{-p}\,du
=c_p^2 |t|^{-2p}\int_{u = 0}^1 u^{-p}(1-u)^{-p}\,du,
\end{equation*}
and the result follows by evaluating the beta integral.
\end{proof}

In particular, recalling $\Gamma(1 / 2) = \sqrt{\pi}$, Lemmas~\ref{L:1/2}
and~\ref{L:p01} yield
\begin{equation}
\label{secbd}
|\phi(t)| \leq \frac{4 \pi}{|t|}.
\end{equation}
This proves~\eqref{phip} for $p = 1$, with $c_1 \leq 4 \pi$, and thus 
by~\eqref{cpmon} for every
$p \le 1$ with $c_p \leq (4 \pi)^p$;
applying \refL{L:p01} again, we obtain the finiteness of~$c_p$
in~\eqref{phip} for all $p < 2$.
Somewhat better numerical bounds are obtained for $1/2 < p < 1$
by taking a geometric average between the cases $p = 1/2$ and $p = 1$:\ 
the inequality
$$
|\phi(t)|\le
(2t^{-1/2})^{2-2p}(4\pi t^{-1})^{2p-1}=2^{2p}\pi^{2p-1}t^{-p},
\qquad t>0,
$$
shows that $c_p \leq 2^{2p} \pi^{2p-1}$, $1/2 \le p \le 1$.
In particular, we have $c_{3/4} \leq \sqrt{8\pi}$, and thus, by
\refL{L:p01},
$c_{3/2} \leq 8\pi^{1/2} \bigl[\Gamma(1/4)\bigr]^2 < 186.4 < 187$.

\begin{lemma}\label{L:p>1}
Let $p > 1$.  Then
$$
c_{p + 1} \le 2^{p+1} c_p^{1 + (1/p)} p / (p - 1).
$$
\end{lemma}

\begin{proof}
Assume that $t \ge 2c_p^{1/p}$. Then,
again using~(\ref{phirecbd}),
\begin{align*}
|\phi(t)|
 &\leq \int^1_{u = 0}\!\min \left( \frac{c_p}{(u t)^{p}}, 1 \right)\,\min
         \left( \frac{c_p}{[(1 - u) t]^{p}}, 1 \right) du \\
 &  =  2 \int^{c_p^{1 /p} t^{-1}}_{u = 0}\!\frac{c_p}{[(1 - u) t]^{p}}\,du + 
\int^{1 - c_p^{1 / p} t^{-1}}_{u = c_p^{1 / p} t^{-1}}\!
\frac{c^2_p}{[u (1 - u) t^2]^{p}}\,du \\
 &\leq \frac{2}{\left[ 1 - c_p^{1 / p} t^{-1} \right]^{p}}
         \frac{c_p^{1 + (1 / p)}}{t^{p +1}} + 
        2 \frac{c^2_p}{t^{2 p}} \int^{1/2}_{u = c_p^{1 /p}t^{-1}}\!
        \frac{du}{[u (1 - u)]^{p}} \\
 &\leq \frac{2}{(1 / 2)^{p}} c_p^{1 + (1 / p)} t^{- (p +1)} 
  +\frac{2}{(1 / 2)^{p}} \frac{c^2_p}{t^{2 p}}
         \int^{1/2}_{u = c_p^{1 / p} t^{-1}}\!u^{- p}\,du \\
 &\leq 2^{p +1} \left\{ c_p^{1 + (1 / p)} t^{- (p +1)}
         +\frac{1}{p - 1} c^2_p t^{- 2 p} \left[ c_p^{1 / p}
         t^{-1} \right]^{- (p - 1)} \right\} \\
 &  =  2^{p +1} c_p^{1 + (1 / p)} \frac{p}{p - 1} t^{- (p +1)} .
\end{align*} 
We have derived the desired bound for all $t \geq 2 c_p^{1 / p}$.
But also, for all $0 < t < 2 c_p^{1 / p}$, we have 
$$
2^{p +1} c_p^{1 + (1 / p)} \frac{p}{p - 1} t^{- (p +1)} 
\ge \frac p{p-1}\ge 1\ge |\phi(t)|,
$$
so the estimate holds for all $t>0$.
\end{proof}

\refL{L:p>1} completes the proof of finiteness of every~$c_p$ in~\eqref{phip}
(by induction),
and of the estimate~\eqref{cp+1}.
The bound for~$c_{3/2}$ obtained above now shows (using {\tt Maple}) that
$c_{5/2} < 103215$, which then gives $c_{7/2} < 197102280$.

We can rewrite \eqref{cp+1} as
\begin{align*}
c_{p+1}^{1/(p+1)}&\le 2c_p^{1/p} \Bigl(1 + \frac{1}{p-1}\Bigr)^{1/(p+1)} 
\le 2 c_p^{1/p} \exp\Bigl(\frac{1}{(p-1)(p+1)}\Bigr)\\
&= 2c_p^{1/p} \exp\Bigl(\frac{1}{2 (p-1)} - \frac{1}{2 (p+1)}\Bigr).
\end{align*}
Hence, by induction, if $p = n+ \frac{5}{2}$ for a nonnegative integer
$n$, then
\begin{equation*}
c_p^{1/p}
\le 2^n c_{5/2}^{2/5} e^{(1/3)+(1/5)} = C 2^p,
\end{equation*}
where $C := 2^{-5/2} e^{8/15} c_{5/2}^{2/5} < 30.6 < 2^5$, using the above estimate
of~$c_{5/2}$.
Consequently, $c_p^{1/p}<2^{p+5}$ when $p=n+\frac{5}{2}$. 
For general $p>\xfrac32$ we now use~\eqref{cpmon} with $p_1=p$ and 
$p_2 = \lceil p -\frac{5}{2} \rceil + \frac{5}{2}$,  obtaining 
$c_p^{1 / p} < 2^{p_2+5}<2^{p + 6}$; 
the case $p \le 3/2$ follows from~\eqref{cpmon} and the
estimate $c_{3/2}^{2/3} < 33 < 2^6$.  
This completes the proof of~\eqref{cpbound} and
hence of \refT{T:phi}.
\end{proof}

\begin{remark}
We used~\eqref{fix} in two different ways.
In the first step we conditioned on the values of~$Y$ and~$Z$, while
in the inductive steps we conditioned on~$U$.
\end{remark}

\begin{remark}
A variety of other bounds are possible.  For example, if we begin with the
inequality~\eqref{secbd}, use~\eqref{phirecbd}, and proceed just as in the proof of
\refL{L:p>1}, we can easily derive the following result in the case $t
\geq 8 \pi$:
\begin{equation}
\label{thirdbd}
|\phi(t)| 
\leq \frac{32 \pi^2}{t^2} \Bigl(\ln \Bigl(\frac{t}{4\pi}\Bigr)+2\Bigr)
\leq \frac{32 \pi^2 \ln t}{t^2}\mbox{\ \ for all $t \geq 1.72$.}
\end{equation}
The result is trivial for $1.72 \leq t < 8 \pi$, since then the bounds
exceed unity. 
\end{remark}

Since~$Y$ has finite moments of all orders [recall~\eqref{finitemgf}],
the characteristic function~$\phi$ is infinitely differentiable.
\refT{T:phi} implies a rapid decrease of all derivatives, too.
\begin{theorem}\label{T:chis}
For each real $p\ge0$ and integer $k\ge0$, there is a constant
$c_{p,k}$ such that
\begin{equation*}
|\phi^{(k)}(t)|\le c_{p,k}|t|^{-p}\mbox{\rm\ \ for all $t \in \RR$.}
\end{equation*}
\end{theorem}
\begin{proof}
The case $k=0$ is \refT{T:phi}, and the case $p = 0$ follows by
$|\phi^{(k)}(t)| \le \EE |Y|^k$. The remaining cases follows from these
cases by induction on~$k$ and the following calculus lemma.
\end{proof}

\begin{lemma}
Suppose that~$g$ is a complex-valued function on $(0,\infty)$ and that
$A,B,p>0$ are such that 
$|g(t)| \le A t^{-p}$ and 
$|g''(t)| \le B$ for all $t>0$.
Then $|g'(t)| \le 2\sqrt{AB} t^{-p/2}$.
\end{lemma}

\begin{proof}
Fix $t>0$ and let $\theta=\arg (g'(t))$. For $s > t$,
\begin{equation*}
\Re(e^{-i \theta} g'(s))
\ge \Re(e^{-i \theta} g'(t)) - |g'(s) - g'(t)|
\ge |g'(t)| - B(s - t)
\end{equation*}
and thus, integrating from $t$ to $t_1 := t + (|g'(t)| / B)$,
\begin{equation*}
\begin{split}
\Re\bigl(e^{-i \theta} (g(t_1) - g(t))\bigr)
&\ge \int_t^{t_1} \bigl(|g'(t)| - B(s - t)\bigr)\,ds\\
&= (t_1 - t) |g'(t)| - \mbox{$\frac{1}{2}$} B (t_1 - t)^2
= |g'(t)|^2 / (2 B).
\end{split}
\end{equation*}
Consequently,
\begin{equation*}
|g'(t)|^2 / (2B) 
\le
|g(t)| + |g(t_1)|
\le2 A t^{-p},
\end{equation*}
and the result follows.
\end{proof}

In other words, the characteristic function~$\phi$ belongs to the
class~$\cS$ of infinitely differentiable functions that, together with
all derivatives, decrease more rapidly than any power.
(This is the important class of test functions for tempered
distributions, introduced by Schwartz~\cite{Schwartz}; it is often
called the class of \emph{rapidly decreasing} $C^\infty$ functions.)

\sect{The limiting Quicksort density~$f$ and its derivatives}
\label{density}

We can now improve the result by
Tan and Hadjicostas~\cite{TanH} 
on existence of a density~$f$ for~$Y$.
It is an immediate consequence of \refT{T:phi}, with $p = 0$ and $p = 2$,
say, that the characteristic function~$\phi$ is integrable over
the real line.  It is well-known---see, 
\eg,\ \cite[Theorem~XV.3.3]{Feller2}---that this implies
that~$Y$ has a bounded continuous density~$f$ given by the
Fourier inversion formula
\begin{equation}
\label{Fourierinv}
f(x) = \frac{1}{2 \pi}\int^{\infty}_{t = -\infty}\!e^{- i tx}\,\phi(t)\,dt,
\ \ x \in \RR.
\end{equation} 
Moreover, using \refT{T:phi} with $p=k + 2$, we see that $t^k \phi(t)$ is also
integrable for each integer $k \ge 0$, which by a standard argument
(\cf{} \cite[Section~XV.4]{Feller2})  
shows 
that~$f$ is infinitely smooth, with a $k$th derivative ($k \geq 0$) given by
\begin{equation}\label{Fourierinvk}
f^{(k)}(x) 
= \frac{1}{2 \pi}\int^{\infty}_{t = -\infty}\!
(- i t)^k\,e^{- i t x}\,\phi(t)\,dt,
\qquad x \in \RR.
\end{equation}
It follows further that the derivatives are bounded, with 
\begin{equation}
\label{fkbd}
\sup_x |f^{(k)}(x)| \le \frac1{2 \pi} \int_{t = -\infty}^\infty |t|^k\,|\phi(t)|\,dt\qquad
\mbox{($k \geq 0$)},
\end{equation}
and these bounds in turn can be estimated using \refT{T:phi}.
Moreover, as is well known~\cite{Schwartz}, 
\cite[Theorem 7.4]{Rudin},
an extension of this argument shows
that the class~$\cS$
discussed at the end of Section~\ref{chf}
is preserved by the Fourier transform, and thus
\refT{T:chis} implies that $f\in\cS$:
\begin{theorem}\label{T:smooth}
The {\tt Quicksort} limiting distribution has an
infinitely differentiable density function~$f$.
For each real $p \ge 0$ and integer $k \ge 0$, there is a constant
$C_{p,k}$ such that
$$
|f^{(k)}(x)|\le C_{p,k}|x|^{-p}\mbox{\rm\ \ for all $x \in \RR$.}
\eqno\qedsymbol
$$
\end{theorem}

For numerical bounds on~$f$, we can use
\eqref{fkbd} with $k = 0$
and \refT{T:phi} for
several different~$p$
(in different intervals);
for example, using $p=0$, $1/2$, $1$, $3/2$, $5/2$, $7/2$,
and taking $t_1 = 4$, $t_2 = 4\pi^2$, $t_3 = (187 / (4 \pi))^2$,
$t_4 = 103215/187$, 
$t_5 = 197102280/103215$, 
\begin{equation}\label{fxest}
\begin{split}
f(x)
&\le
\frac1{2\pi}\int_{t = -\infty}^\infty
|\phi(t)|\,dt
=\frac1{\pi}\int_{t = 0}^\infty
|\phi(t)|\,dt\\
&\le
\frac1{\pi} 
\int_{t = 0}^{\infty}
\min(1,2t^{-1/2},4\pi t^{-1},187 t^{-3/2}, 
103215 \,t^{-5/2}, 
197102280\, t^{-7/2})\,dt \\
&=\frac1{\pi}
\biggl(
\int_{t = 0}^{t_1} dt
+\int_{t = t_1}^{t_2}
2t^{-1/2}\,dt
+\int_{t = t_2}^{t_3}
4\pi t^{-1}\,dt
+\int_{t = t_3}^{t_4}
187\, t^{-3/2}\,dt\\
&\qquad\qquad\qquad+\int_{t = t_4}^{t_5}
103215 \,t^{-5/2}\,dt 
+\int_{t = t_5}^{\infty}
197102280\, t^{-7/2}\,dt 
\biggr)\\
&\le 18.2.
\end{split}
\end{equation}

\begin{remark}
\label{R:trick}
We can do somewhat better by using  the first bound in~\eqref{thirdbd}
over the interval $(103.18, 1984)$
instead of  (as above) \refT{T:phi} with 
$p=1$, $3 / 2$, $5 / 2$,  $7 / 2$ over $(103.18,t_3)$, 
$(t_3, t_4)$, 
$(t_4, t_5)$, 
$(t_5, 1984)$, respectively.
This gives
$$
f(x) < 15.3.
$$ 
\end{remark}

Similarly,
\eqref{fkbd} with $k = 1$ and the same
estimates of $|\phi(t)|$ as in~\eqref{fxest} yield
$$
|f'(x)|
\le
\frac1{2\pi}\int_{t = -\infty}^\infty
|t| |\phi(t)|\,dt
=\frac1{\pi}\int_{t = 0}^\infty
t|\phi(t)|\,dt
< 3652.1,
$$
which can be reduced to $2492.1$ by proceeding as in \refR{R:trick}.
The bound can be further improved to~$2465.9$ by using also $p =
9/2$.

Somewhat better bounds are obtained by using more values of~$p$ in the
estimates of the integrals, but the improvements obtained in this way
seem to be slight. 
We summarize the bounds we have obtained.
\begin{theorem}\label{T:maxf}
The limiting {\tt Quicksort} density function~$f$ satisfies
$\max_x f(x) < 16$ and $\max_x |f'(x) | < 2466$.
\nopf
\end{theorem}

The numerical bounds obtained here are far from sharp;
examination of Figure~4 of~\cite{TanH} suggests that $\max f < 1$ and $\max |f'| < 2$.
Our present technique cannot hope to produce a better bound on~$f$ than~$4 / \pi > 1.27$, since
neither \refL{L:hyz} nor~\eqref{phirecbd} can improve on the bound $|\phi(t)| \leq 1$
for $|t| \leq 4$.  Further, \emph{no} technique based on~\eqref{fkbd} can hope to do
better than the actual value of $(2 \pi)^{-1} \int^{\infty}_{t = -\infty}\!|\phi(t)|\,dt$,
which from cursory examination of Figure~6 of~\cite{TanH} appears to be about~$2$.

\section{An integral equation for the density~$f$} \label{finteq}

Our estimates are readily used to justify rigorously the following functional equation.

\begt\label{T:inteq}
The continuous limiting {\tt Quicksort} density~$f$ satisfies
\rompar(pointwise) the integral equation
$$
f(x) = \int^1_{u = 0} \int_{y \in \RR} f(y)\,f \left( \frac{x - g(u) -
(1 - u) y}{u} \right) 
\frac{1}{u}\,dy\,du,\qquad x \in \RR,
$$
where~$g(\cdot)$ is as in~\eqref{gu}.
\ent
\begin{proof}
For each~$u$ with $0<u<1$,
the random variable
\begin{equation}
\label{fixurv}
u Y + (1 - u) Z + g(u)
\end{equation}
[with notation as in~(\ref{fix})] 
has the density function
\begin{equation}
\label{dens2}
f_u(x):= \int_{z \in \RR} f(z)\,f \left( \frac{x - g(u) - (1 - u)
z}{u} \right) 
\frac{1}{u}\,dz,
\end{equation}
where the integral converges for each~$x$ since,
using \refT{T:maxf},
the integrand is bounded by
$f(z) (\max f)/ u \le 16 f(z) / u$; 
dominated convergence using the continuity of~$f$ and
the same bound 
shows further that~$f_u$ is continuous.

This argument yields the bound $f_u(x) \le 16 / u$, and since 
$f_u=f_{1-u}$ by symmetry in~\eqref{fixurv}, we have
$f_u(x) \le 16/ \max(u, 1 - u) \le 32$.  This uniform bound,
\eqref{fix}, and dominated convergence again imply that
$\int_0^1 f_u(x)\,du$ is a continuous density for~$Y$,
and thus equals $f(x)$ for every~$x$.
\end{proof}

It was shown in~\cite{TanH} that~$f$ is positive almost everywhere;
we now can improve this by removing the qualifier ``almost.''

\begin{corollary}
\label{C:positive} 
The continuous limiting {\tt Quicksort} density function is
everywhere positive.
\end{corollary}
\begin{proof}
We again use the notation~\eqref{dens2} from the proof of \refT{T:inteq}.
Fix $x \in \RR$ and $u \in (0, 1)$.
Since~$f$ is almost everywhere positive~\cite{TanH},
the integrand in~\eqref{dens2} is positive almost everywhere.
Therefore $f_u(x) > 0$.  Now we integrate over $u \in (0, 1)$
to conclude that $f(x) > 0$.
\end{proof}

Alternatively, \refC{C:positive} can be derived directly from \refT{T:inteq},
without recourse to~\cite{TanH}.  Indeed, if $f(y_0) > 0$ and $u_0 \in (0, 1)$,
set $x = y_0 + g(y_0)$; then the integrand in the double integral for~$f(x)$ in
\refT{T:inteq} is postive for $(u, y)$ equal to $(u_0, y_0)$, and therefore, by
continuity, also in some small neighborhood thereof.  It follows that
$f(y_0 + g(u_0)) > 0$.  Since $u_0$ is arbitrary and the image of~$(0, 1)$
under~$g$ is $( - (2 \ln 2 - 1), 1)$, an open interval containing the origin,
\refC{C:positive} follows readily.    

\begin{remark}
In future work, we will use arguments similar to those of this paper,
together with other arguments, to show that when one applies the method
of successive substitutions to the integral equation in \refT{T:inteq},
the iterates enjoy exponential-rate uniform convergence to~$f$.  This will
settle an issue raised in the third paragraph of Section~3 in~\cite{TanH}.
\end{remark}

\end{document}